\documentclass[11pt]{article}
\usepackage[a4paper,margin=25mm]{geometry}
\usepackage{amsmath,amsthm}
\usepackage{enumitem,microtype}
\usepackage{booktabs,multirow,hhline,makecell}
\usepackage{hyperref,xcolor}
\usepackage{tikz}
\usepackage[nameinlink,noabbrev]{cleveref}

\newtheorem{theorem}{Theorem}[section]
\newtheorem{proposition}[theorem]{Proposition}
\newtheorem{lemma}[theorem]{Lemma}
\newtheorem{claim}{Claim}
\theoremstyle{definition}
\newtheorem{definition}[theorem]{Definition}
\newcommand{\cE}{\mathcal E}
\title{Minimum degree thresholds for even cycle residues}

\author{Hojin Chu
\thanks{Korea Institute for Advanced Study, Seoul 02455, Republic of Korea. Email: \texttt{hojinchu@kias.re.kr}. Supported by a KIAS Individual Grant (CG101801) at Korea Institute for Advanced Study.}
\and Masaki Kashima
\thanks{Faculty of Science and Technology, Keio University. Email: \texttt{masaki.kashima10@gmail.com}. Supported by JSPS KAKENHI Grant Number JP25KJ2077.}
\and Shun-ichi Maezawa\thanks{Department of Information Science, College of Humanities and Sciences, Nihon University. Email: \texttt{maezawa.mw@gmail.com}. Supported by JSPS KAKENHI Grant Number JP25K17301.}
\and Boram Park
\thanks{Department of Mathematics Education, Seoul National University, Seoul 08826, Republic of Korea. Email: \texttt{borampark@snu.ac.kr}. Supported by the National Research Foundation of Korea (NRF) grant funded by the Korea government (MSIT) (No.~RS-2025-00523206) and by the New Faculty Startup Fund from Seoul National University.}
}
\date{}

\begin{document}

\maketitle
\begin{abstract}
Dean's conjecture asserts that, for every integer $m\ge3$, every graph of minimum degree at least $m$ has a cycle whose length is divisible by $m$. This degree bound is best possible for odd $m$. For every integer $k\ge2$, we prove that every graph of minimum degree at least $2k-1$ has a cycle whose length is congruent to $r$ modulo $2k$ for every even $r\not\equiv2\pmod{2k}$. Consequently, $2k-1$ is the sharp minimum degree threshold for a cycle whose length is divisible by $2k$. For bipartite graphs, we prove that minimum degree at least $k$ guarantees the same residue classes for every $k\ge4$, and we characterize their end-blocks when no cycle has length congruent to $2$ modulo $2k$. The bipartite degree bound is also sharp, and the exclusion of residue $2$ is necessary even for connected graphs of arbitrarily large order.
\end{abstract}



\section{Introduction}

All graphs considered in this paper are finite and simple. For integers $m\ge1$ and $\ell$, an \emph{$(\ell\bmod m)$-cycle} is a cycle whose length is congruent to $\ell$ modulo $m$.

Dean's conjecture asserts that, for every integer $m\ge3$, every graph of minimum degree at least $m$ has a $(0\bmod m)$-cycle. The conjecture was proved for $m=3$ by Chen and Saito~\cite{CS}, for $m=4$ by Dean, Lesniak, and Saito~\cite{DLS}, and for every $m\ge6$ by Luo, Ma, and Zhao~\cite{LMZ}. A proof of the remaining case $m=5$ was recently announced by Sieniawski~\cite{Sieniawski}. For odd $m$, the graph $K_{m-1,N}$ with $N\ge m-1$ has minimum degree $m-1$ and no $(0\bmod m)$-cycle, so the degree bound cannot be reduced.

Thomassen~\cite{Thomassen} conjectured that every graph of minimum degree at least $m+1$ has an $(r\bmod m)$-cycle for every even integer $r$. Liu and Ma~\cite{LM} proved this for even $m$, and Gao, Huo, Liu, and Ma~\cite{GHLM} proved it for all $m\ge3$. Luo, Ma, and Zhao~\cite{LMZ} showed that, for $m\ge6$, minimum degree at least $m$ guarantees every even residue class other than $2$. For even moduli $2k\ge4$, we lower the degree bound to $2k-1$ and prove a corresponding result for bipartite graphs. \Cref{tab:thresholds} compares our results with previous ones.

We first consider bipartite graphs. Liu and Ma~\cite{LM} proved that a bipartite graph of order at least two in which all but at most one vertex have degree at least $k$ has $k-1$ cycles of consecutive even lengths. Bai, Grzesik, Li, and Prorok~\cite{BGLP} showed that every $2$-connected bipartite graph of minimum degree at least $k\ge4$ has $k$ such cycles, unless it is isomorphic to $K_{k,n}$ for some $n\ge k$. Without the connectivity assumption, Luo, Ma, and Zhao~\cite{LMZ} showed for $k\ge7$ that every graph of minimum degree at least $k$ has $k$ cycles whose lengths form an arithmetic progression with common difference one or two, unless every end-block belongs to an explicit family of graphs.

A block $B$ of a graph $G$ is an \emph{end-block} if $B$ has at most one cut-vertex of $G$. A connected component that is $2$-connected is also an end-block. For integers $2\le t\le k<n$, let $H_{k,n;t}$ be obtained from $K_{k,n}$ by deleting $k-t$ edges incident with one vertex in the part of size $n$, and let
\[\cE_k=\{K_{k,k}\}\cup\{H_{k,n;t}:2\le t\le k<n\}.\]
Every member of $\cE_k$ has cycles of exactly the lengths $4,6,\ldots,2k$. For bipartite graphs, the exceptional family of Luo, Ma, and Zhao is exactly $\cE_k$. Our first result extends the corresponding conclusion for cycle lengths modulo $2k$ from $k\ge7$ to every $k\ge4$.

\begin{theorem}\label{thm:main}
Let $k\ge4$ be an integer, and let $G$ be a bipartite graph with $\delta(G)\ge k$. Then exactly one of the following holds.
\begin{enumerate}[label=\textnormal{(\roman*)},leftmargin=9mm]
\item $G$ has an $(r\bmod 2k)$-cycle for every even integer $r$.
\item Every end-block of $G$ is isomorphic to a member of $\cE_k$, and every non-end-block of $G$ has no $(2\bmod 2k)$-cycle.
\end{enumerate}
\end{theorem}

Thus $2$ is the only even residue class modulo $2k$ that can be absent in a bipartite graph of minimum degree at least $k$. The bipartite graph $K_{k-1,N}$ with $N\ge k$ has minimum degree $k-1$ and no $(0\bmod 2k)$-cycle, showing that this degree bound is sharp. Alternative~(ii) is also necessary, since identifying one vertex from each of arbitrarily many copies of $K_{k,k}$ gives connected bipartite graphs of arbitrarily large order and minimum degree $k$ with no $(2\bmod 2k)$-cycle. Increasing the minimum degree to $k+1$ guarantees $k$ cycles of consecutive even lengths~\cite{GHLM}, and hence every even residue class modulo $2k$.

For arbitrary graphs, Bai, Grzesik, Li, and Prorok~\cite{BGLP} proved that, for every integer $k\ge2$, a $2$-connected graph of order at least $2k+2$ and minimum degree at least $2k-1$ has an $(r\bmod 2k)$-cycle for every even integer $r$. Our second result removes the connectivity and order assumptions for all even residues other than $2$.

\begin{theorem}\label{thm:general}
Let $k\ge2$ be an integer. Every graph $G$ with $\delta(G)\ge2k-1$ has an $(r\bmod 2k)$-cycle for every even integer $r\not\equiv2\pmod{2k}$.
\end{theorem}

For $k=2$, \Cref{thm:general} is the theorem of Dean, Lesniak, and Saito~\cite{DLS}. For $k\ge3$, it lowers the degree bound $2k$ of Luo, Ma, and Zhao~\cite{LMZ} by one. Chen, Liu, Wang, Wei, and Yang~\cite{CLWWY} obtained a result on consecutive even cycle lengths in terms of the number of edges, which also implies \Cref{thm:general} for sufficiently large $k$. The degree bound $2k-1$ is sharp, since $K_{2k-1}$ has minimum degree $2k-2$ and no $(0\bmod 2k)$-cycle. Residue $2$ behaves differently: identifying one vertex from each of arbitrarily many copies of $K_{2k+1}$ gives connected graphs of arbitrarily large order and minimum degree $2k$ in which every cycle has length at most $2k+1$, so they have no $(2\bmod 2k)$-cycle. Hence, for residue $2$, the bound $2k+1$ of~\cite{GHLM} is best possible.

\begin{table}[htbp]
\centering
\footnotesize
\renewcommand{\cellalign}{l}
\setcellgapes{3pt}\makegapedcells
\setlength{\aboverulesep}{0pt}
\setlength{\belowrulesep}{0pt}
\setlength{\doublerulesep}{1.3pt}
\setlength{\cmidrulewidth}{0.2pt}\setlength{\arrayrulewidth}{0.3pt}
\begin{tabular}{l|l|l|l|l}
\toprule
\multirow{2}{*}{} & \multirow{2}{*}{\makecell{$(r\bmod 2k)$-cycles}} & \multicolumn{2}{c|}{Bound on $\delta(G)$} & \multirow{2}{*}{Sharpness example} \\
\cmidrule{3-4}
 & & Previous result & This paper & \\
\hhline{=|=|=|=|=}
\makecell{General\\[3pt] graphs} & all even $r\not\equiv2$ & \makecell{$2k$ if $k\ge3$~\cite{LMZ}\\[5pt] $3$ if $k=2$~\cite{DLS}} & $2k-1$ if $k\ge2$ (Thm.~\ref{thm:general}) & $K_{2k-1}$ \\
\midrule
\multirow{2}{*}{\makecell{Bipartite\\[3pt] graphs}} & all even $r\not\equiv2$ & \makecell{$k+1$ if $k\ge3$~\cite{GHLM}\\[5pt] $k$ if $k\ge7$~\cite{LMZ}\\[5pt] $k$ if $k\ge4$ and $2$-conn.~\cite{BGLP}} & \multirow{2}{*}[-4ex]{\makecell{$k$ if $k\ge4$ (Thm.~\ref{thm:main});\\[5pt] for all even $r$, unless\\[2pt] \quad every end-block is in $\cE_k$}} & $K_{k-1,N}$  ($N\ge k$) \\
\cmidrule{2-3}\cmidrule{5-5}
 & all even $r$ & \makecell{$k+1$ if $k\ge3$~\cite{GHLM}\\[5pt] $k$ if $k\ge7$, unless every\\[2pt] \quad end-block is in $\cE_k$~\cite{LMZ}} & &  $K_{k,k}$ \\
\bottomrule
\end{tabular}
\caption{Minimum degree conditions guaranteeing an $(r\bmod 2k)$-cycle for every $r$ in the indicated range.}\label{tab:thresholds}
\end{table}

The paper is organized as follows. In Section~\ref{sec:preliminaries}, we introduce notation and collect preliminary results. In Section~\ref{sec:admissible}, we prove the main lemmas on admissible cycles in $k$-weak bipartite graphs. In Section~\ref{sec:main-proofs}, we use these lemmas to prove \Cref{thm:main,thm:general}.

\section{Preliminaries}\label{sec:preliminaries}

For a graph $G$, we write $V(G)$ and $E(G)$ for its vertex and edge sets, respectively. For $v\in V(G)$, let $N_G(v)=\{u\in V(G):uv\in E(G)\}$ and $d_G(v)=|N_G(v)|$ denote its neighborhood and degree, respectively. The minimum degree of $G$ is denoted by $\delta(G)$. For $X\subseteq V(G)$, let $N_G(X)$ be the set of vertices in $V(G)\setminus X$ having a neighbor in $X$.

For distinct vertices $x,y\in V(G)$, let $G+xy$ denote the graph obtained from $G$ by adding the edge $xy$ if it is absent.

For disjoint vertex sets $X$ and $Y$ of a graph $G$, an \emph{$(X,Y)$-path} is a path having one endpoint in $X$, the other endpoint in $Y$, and no internal vertex in $X\cup Y$. We abbreviate $(\{x\},Y)$ to $(x,Y)$ and $(\{x\},\{y\})$ to $(x,y)$. Paths or cycles are called \emph{admissible} if their lengths form an arithmetic progression with common difference one or two, and for paths, we also require the shortest length to be at least two. In a bipartite graph, the common difference is two for admissible cycles and for admissible $(x,y)$-paths with $x,y$ fixed.

\begin{theorem}[Chiba--Ota--Yamashita~\cite{COY}]\label{thm:paths}
Let $m$ be a positive integer, let $x,y$ be distinct vertices of a graph $G$ such that $G+xy$ is $2$-connected, and let $z$ be a vertex of $G$ (possibly $z\in\{x,y\}$) such that $V(G)\setminus\{x,y,z\}\ne\emptyset$. Suppose that $d_G(v)\ge m+1$ for every $v\in V(G)\setminus\{x,y,z\}$. Then $G$ has $m$ admissible $(x,y)$-paths.
\end{theorem}

In applications, if the indicated exceptional vertex is absent, we take $z$ to be one of the two prescribed endpoints.

The following definition is from \cite{LMZ}, where $\delta_2(G)$ is the second smallest degree in a graph $G$.

\begin{definition}
Let $k\ge3$ be an integer. A graph $G$ is \emph{$k$-weak} if one of the following holds.
\begin{enumerate}[label=\textnormal{Type \Roman*.},leftmargin=22mm]
\item $G$ is $3$-connected and $\delta_2(G)\ge k$.
\item $G$ has a unique vertex $\theta$ of degree two, with neighbors $x,y$, such that $G-\theta+xy$ is $3$-connected and has minimum degree at least $k$.
\end{enumerate}
\end{definition}

\begin{lemma}\label{lem:weak}
Let $k\ge3$, let $G$ be a $2$-connected bipartite graph, and let $c\in V(G)$ satisfy $d_G(v)\ge k$ for every $v\in V(G)\setminus\{c\}$. Then either $G$ is $k$-weak or $G$ has $k$ admissible cycles.
\end{lemma}

\begin{proof}
Suppose that $G$ contains no $k$ admissible cycles. If $G$ is $3$-connected, then it is $k$-weak of Type~I. Otherwise, $G$ has a $2$-cut $\{x,y\}$ for some $x,y\in V(G)$.

Suppose that the components of $G-\{x,y\}$ can be partitioned into two unions $U_1, U_2$, each of order at least two. For $i=1,2$, let $G_i=G[U_i\cup\{x,y\}]$. Since $G$ is $2$-connected, $G_i+xy$ is $2$-connected. Every vertex of $U_i\setminus\{c\}$ has degree at least $k$ in $G_i$. By \Cref{thm:paths}, each $G_i$ has $k-1$ admissible $(x,y)$-paths. Joining one path from each side gives $2k-3\ge k$ admissible cycles, a contradiction.

A singleton component of $G-\{x,y\}$ has degree at most two in $G$, and hence must be $\{c\}$. It follows that $d_G(c)=2$, $N_G(c)=\{x,y\}$, and $G-\{x,y\}$ has precisely two components, one of which is $\{c\}$. Since $G$ is bipartite, $xy\notin E(G)$. Then the graph $H:=G-c+xy$ is $2$-connected and satisfies  $\delta(H)\ge k$. Note that $N_G(c)$ is the unique $2$-cut of $G$ since $\{x,y\}$ was arbitrarily chosen. Since $H-\{x,y\}$ is connected, $H$ is $3$-connected. Therefore $G$ is $k$-weak of Type~II.
\end{proof}

We also use the following results.

\begin{theorem}[Lu--Yu~\cite{LY}]\label{thm:LY}
Let $G$ be a graph with $\delta(G)\ge3$. If $G$ has no $(1\bmod 3)$-cycle, then $G$ contains the Petersen graph as an induced subgraph.
\end{theorem}

\begin{theorem}[Chen--Saito~\cite{CS}]\label{thm:CS}
Let $G$ be a graph with at least two vertices. If all but at most one vertex of $G$ have degree at least three, then $G$ has a $(0\bmod 3)$-cycle.
\end{theorem}

For a graph $G$ and a block $B$ of $G$, a vertex of $B$ that is not a cut-vertex of $G$ is called an \emph{inner vertex} of $B$ in $G$.

\begin{lemma}[Bai--Grzesik--Li--Prorok~\cite{BGLP}]\label{lem:auxiliary}
Let $x,y$ be distinct vertices of a graph $G$ such that $Q=G-\{x,y\}$ is connected and every end-block $B$ of $Q$ has an inner vertex adjacent to $x$ or $y$. If there exist two vertex-disjoint edges from $\{x,y\}$ to $V(G)\setminus\{x,y\}$, then $G+xy$ is $2$-connected.
\end{lemma}

\begin{lemma}[Bai--Grzesik--Li--Prorok~\cite{BGLP}]\label{lem:contraction}
Let $G$ be a $3$-connected graph, let $C$ be a cycle of $G$, and let $Q$ be a component of $G-V(C)$. The graph obtained from $G$ by contracting $Q$ to a single vertex is $3$-connected.
\end{lemma}

\section{Admissible cycles in $k$-weak bipartite graphs}\label{sec:admissible}

Throughout the section, let $G$ be a $2$-connected bipartite $k$-weak graph for some $k\ge3$, and fix $c\in V(G)$ such that $d_G(v)\ge k$ for every $v\in V(G)\setminus\{c\}$. If $G$ has a $2$-cut, then $G$ is of Type~II, and its only $2$-cut is $N_G(c)$, whose deletion leaves $\{c\}$ and one other component.

\subsection{The case with a cycle of length four}

For $n\ge2$, let $K_{n,n}^{-}$ denote the graph obtained from $K_{n,n}$ by deleting one edge.

\begin{lemma}\label{lem:four-cycle}
If $G$ has a cycle of length four, then $G$ has $k$ admissible cycles or $G\in\cE_k$.
\end{lemma}

\begin{proof}
Suppose that $G\notin\cE_k$ and that $G$ contains no $k$ admissible cycles. Choose $s$ as large as possible such that $G$ has a subgraph isomorphic to $K_{s,s}$. If $G$ has a subgraph isomorphic to $K_{s+1,s+1}^{-}$, let $F$ be such a subgraph. Otherwise, choose a subgraph $F$ of $G$ isomorphic to $K_{s,t}$ with $t\ge s$ as large as possible. Let $(S,T)$ be the bipartition of $F$, and let $r=|S|\le |T|$. Since $G$ has a cycle of length four, $r\ge 2$. Note that $r\in\{s,s+1\}$.

\begin{claim}\label{cl:core-neighbors}
For every $v\in V(G)\setminus V(F)$, $|N_G(v)\cap V(F)|\le r-1$.
\end{claim}

\begin{proof}
Suppose to the contrary that $|N_G(v)\cap V(F)|\ge r$ for some $v\in V(G)\setminus V(F)$. Since $G$ is bipartite, $v$ cannot have neighbors in both $S$ and $T$. If $F\cong K_{s+1,s+1}^{-}$, then together with $v$, it  contains $K_{s+1,s+1}$ as a subgraph, contrary to the choice of $s$. Suppose $F\cong K_{s,t}$. If $t=s$, then $F$ and $v$ form $K_{s,s+1}$, contrary to the choice of $t$. Thus $t>s$. If $v$ has neighbors in $S$, then it contradicts the choice of $t$. If $v$ has neighbors in $T$, then together with $v$, it contains $K_{s+1,s+1}^{-}$ as a subgraph, a contradiction.
\end{proof}

If $r>k$, then $F$ contains cycles of lengths $4,6,\ldots,2k+2$, which are $k$ admissible cycles. Thus $r\le k$. By Claim~\ref{cl:core-neighbors}, a singleton component of $G-V(F)$ must be $\{c\}$.

We show that $G-V(F)$ has a component $Q$ which is not a singleton.
Suppose to the contrary that $V(G)\setminus V(F)\subseteq\{c\}$. If $F\cong K_{r,r}^-$, one endpoint of its missing edge, different from $c$, has degree $r-1<k$ in $G$, a contradiction. Thus $F\cong K_{r,t}$, with $|S|=r\le t=|T|$. Both parts of $G$ have order at least $k$ by the degree condition, since each part contains a vertex different from $c$. If one part of $G$ has order $k$, every vertex of the other part except possibly $c$ is adjacent to all of it, which implies that $G\in\cE_k$. Thus both parts of $G$ have order at least $k+1$. As $|S|=r\le k$, this forces $r=k$ and $G-c=F\cong K_{k,t}$ with $t>k$, where $c$ is contained in the part containing $S$. Since $d_G(c)\ge2$, $G$ has cycles of lengths $4,6,\ldots,2k,2k+2$, so $G$ has $k$ admissible cycles, a contradiction. Hence $G-V(F)$ has a component $Q$ which is not a singleton.

\vspace{3pt}
\noindent\textbf{Case 1.} $Q$ has neighbors in both $S$ and $T$.

If $F$ is a complete bipartite graph, let $T'=T$. Otherwise, let $s_0t_0$ be the missing edge, where $s_0\in S$ and $t_0\in T$. Since $G$ is $k$-weak and $Q$ is not a singleton, $Q$ has a neighbor in $V(F)\setminus\{s_0,t_0\}$. After interchanging $S$ and $T$ if necessary, we may therefore choose $T'=T\setminus\{t_0\}$ and assume $N_G(V(Q))\cap T'\ne\emptyset$. For every $s\in S$ and $t\in T'$, $F$ has $(s,t)$-paths of lengths $1,3,\ldots,2r-1$, where the longest one passes through $t_0$ if $F\cong K_{r,r}^-$, as $r\ge3$ in this case. It suffices to find $k-r+1$ admissible $(S,T')$-paths with all internal vertices in $Q$, since together with the $(s,t)$-paths in $F$, they yield $k$ admissible cycles of $G$, a contradiction.

We construct a graph $H$ from $Q$ by adding two nonadjacent vertices $s^*$ and $t^*$, adjacent to the vertices in $N_G(S)\cap V(Q)$ and $N_G(T')\cap V(Q)$, respectively. 
If $F\cong K_{r,r}^{-}$ and $c$ is a leaf of $Q$ whose only neighbor in $F$ is $t_0$, delete $c$ from $H$ and let $z$ be its neighbor in $Q$; otherwise, let $z=c$; see Figure~\ref{fig:auxiliary-graphs}(a). Then $H$ is bipartite with $s^*$ and $t^*$ in different parts, $H-\{s^*,t^*\}$ is connected, and both $s^*$ and $t^*$ have neighbors in $H$. Since the neighbors of $s^*$ and $t^*$ lie in different parts, $H-\{s^*,t^*\}$ has at least two vertices, so $|V(H)|\ge4$.

Suppose that some end-block $D$ of $H-\{s^*,t^*\}$ has no inner vertex adjacent to $s^*$ or $t^*$. Since both $s^*$ and $t^*$ have neighbors in $H$, $D$ has a cut-vertex $b$ of $H$. If $F$ is complete bipartite, then $b$ is a cut-vertex of $G$, a contradiction. Otherwise, $\{b,t_0\}$ separates the inner vertices of $D$ from $S$ in $G$. By $k$-weakness, $D$ is the edge $bc$ with $N_G(c)=\{b,t_0\}$, contrary to the construction of $H$. Thus every end-block of $H-\{s^*,t^*\}$ has an inner vertex adjacent to $s^*$ or $t^*$.   By \Cref{lem:auxiliary}, $H+s^*t^*$ is $2$-connected. By Claim~\ref{cl:core-neighbors} and the construction of $H$, $d_H(v)\ge k-r+2$ for every $v\in V(H)\setminus\{s^*,t^*,z\}$. By \Cref{thm:paths}, $H$ has $k-r+1$ admissible $(s^*,t^*)$-paths, which correspond to $k-r+1$ admissible $(S,T')$-paths with all internal vertices in $Q$.

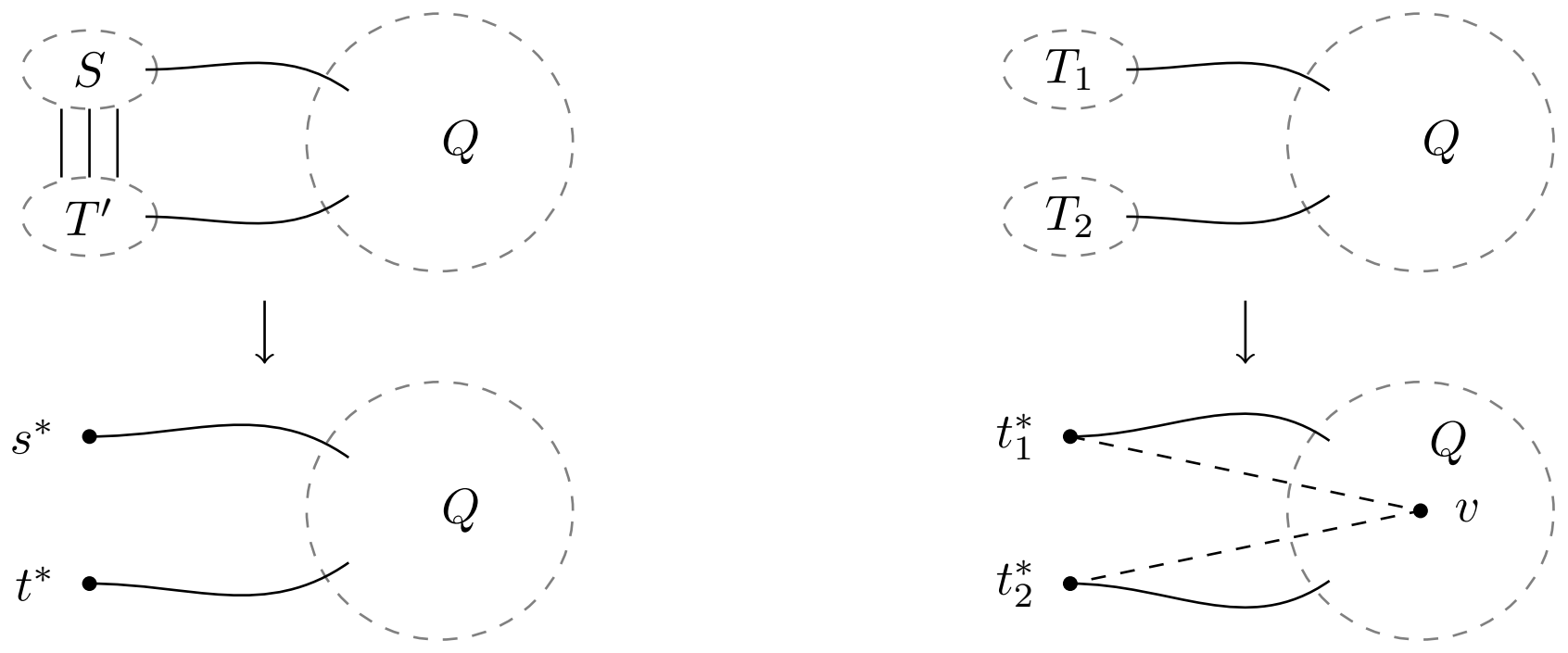
\begin{figure}[htbp]
\centering
\begin{tikzpicture}[font=\small, line width=.5pt]
\begin{scope}
\draw[dashed, gray] (0,3.7) ellipse (.48 and .28);
\draw[dashed, gray] (0,2.65) ellipse (.48 and .28);
\node at (0,3.7) {$S$};
\node at (0,2.65) {$T'$};
\foreach \x in {-.2,0,.2} {\draw (\x,3.42)--(\x,2.93);}
\draw[dashed, gray] (2.5,3.18) ellipse (.95 and .92);
\node at (2.65,3.18) {$Q$};
\draw (.4,3.7) to[out=0,in=145] (1.85,3.55);
\draw (.4,2.65) to[out=0,in=-145] (1.85,2.8);
\draw[->] (1.25,2.05)--(1.25,1.6);
\draw[dashed, gray] (2.5,.55) ellipse (.95 and .92);
\node at (2.65,.55) {$Q$};
\fill (0,1.08) circle (1.5pt) node[left=3pt] {$s^*$};
\fill (0,.03) circle (1.5pt) node[left=3pt] {$t^*$};
\draw (0,1.08) to[out=0,in=145] (1.85,.93);
\draw (0,.03) to[out=0,in=-145] (1.85,.18);
\node at (1.5,-.8) {(a) Case 1};
\end{scope}
\begin{scope}[xshift=7cm]
\draw[dashed, gray] (0,3.7) ellipse (.48 and .28);
\draw[dashed, gray] (0,2.65) ellipse (.48 and .28);
\node at (0,3.7) {$T_1$};
\node at (0,2.65) {$T_2$};
\draw[dashed, gray] (2.5,3.18) ellipse (.95 and .92);
\node at (2.65,3.18) {$Q$};
\draw (.4,3.7) to[out=0,in=145] (1.85,3.55);
\draw (.4,2.65) to[out=0,in=-145] (1.85,2.8);
\draw[->] (1.25,2.05)--(1.25,1.6);
\draw[dashed, gray] (2.5,.55) ellipse (.95 and .92);
\node at (2.7,1.03) {$Q$};
\fill (0,1.08) circle (1.5pt) node[left=3pt] {$t_1^*$};
\fill (0,.03) circle (1.5pt) node[left=3pt] {$t_2^*$};
\draw (0,1.08) to[out=0,in=145] (1.85,1.05);
\draw (0,.03) to[out=0,in=-145] (1.85,.05);
\fill (2.5,.55) circle (1.5pt) node[right=3pt] {$v$};
\draw[dashed] (0,1.08)--(2.5,.55)--(0,.03);
\node at (1.5,-.8) {(b) Case 2};
\end{scope}
\end{tikzpicture}
\caption{Construction of $H$ from $Q$ in Case~1 (where $c$ may also be deleted) and in Case~2.}\label{fig:auxiliary-graphs}
\end{figure}

\noindent\textbf{Case 2.} $Q$ has neighbors in only one of $S$ and $T$.

Let $U\in\{S,T\}$ be such that $N_G(V(Q))\cap V(F)\subseteq U$. Since $G$ is $2$-connected, we can choose a partition $U=T_1\cup T_2$ such that two vertex-disjoint edges join $Q$ to $T_1$ and $T_2$, respectively. Moreover, every end-block of $Q$ has an inner vertex with a neighbor in $U$. We construct a graph $H$ from $Q$ by adding two vertices $t_1^*,t_2^*$, where $v\in V(Q)$ is adjacent to $t_i^*$ if $N_G(v)\cap T_i\ne\emptyset$ or $|N_G(v)\cap U|\ge2$; see Figure~\ref{fig:auxiliary-graphs}(b). Then, by \Cref{lem:auxiliary}, $H+t_1^*t_2^*$ is $2$-connected. Moreover, by Claim~\ref{cl:core-neighbors}, $d_H(v)\ge\min\{k,k-r+3\}$ for every $v\in V(Q)\setminus\{c\}$.

By construction, every $(t_1^*,t_2^*)$-path in $H$ can be replaced by a path of the same length in $G$ with two distinct endpoints in $U$ and all internal vertices in $Q$.
If $r=2$, then $F$ is complete bipartite and $|U|\ge3$, since otherwise $U$ would be a $2$-cut separating $Q$ from $V(F)\setminus U$, contradicting the fact that $G$ is $k$-weak (recall $|V(Q)|\ge2$). By \Cref{thm:paths}, $H$ has $k-1$ admissible $(t_1^*,t_2^*)$-paths. The corresponding paths in $G$, combined with paths of lengths $2,4$ in $F$, form $k$ admissible cycles, a contradiction.

Suppose $r\ge3$. By \Cref{thm:paths}, $H$ has $k-r+2$ admissible $(t_1^*,t_2^*)$-paths. For every distinct $x,y\in U$, $F$ has $(x,y)$-paths of lengths $2,4,\ldots,2r-2$, whether $U=S$ or $U=T$. Combining these with the corresponding paths in $G$ forms $(k-r+2)+(r-1)-1=k$ admissible cycles of $G$, a contradiction.
\end{proof}

\subsection{The case without a cycle of length four}

\begin{lemma}\label{lem:no-four-cycle}
Let $k\ge4$. If $G$ has no cycle of length four, then $G$ has $k$ admissible cycles.
\end{lemma}

\begin{proof}
Suppose to the contrary that $G$ has no $k$ admissible cycles. Let $C$ be a shortest cycle of $G$. If there is one containing $c$, we choose $C$ so that $c\in V(C)$.
Write $|C|=2s\ge6$, and let $W=G-V(C)$. By the choice of $C$, it is induced, and every vertex of $W$ has at most one neighbor on $C$. Thus $d_W(v)\ge k-1$ for every $v\in V(W)\setminus\{c\}$. If $W$ were empty, every vertex of $G$ would have degree two, a contradiction. Moreover, $W$ has no singleton component.

\begin{claim}\label{cl:component-paths}
For every component $Q$ of $W$, there exists $u\in V(Q)$ with a neighbor on $C$ such that, for every $v\in V(Q)\setminus\{u\}$, $Q$ has $k-2$ admissible $(u,v)$-paths.
\end{claim}

\begin{proof}
If $Q$ is $2$-connected, choose $u\in V(Q)$ with a neighbor on $C$. By \Cref{thm:paths}, with exceptional vertex $c$, $Q$ has $k-2$ admissible $(u,v)$-paths for every $v\in V(Q)\setminus\{u\}$.

Suppose that $Q$ is not $2$-connected. Choose an end-block $E$ of $Q$ such that $c$ is not an inner vertex of $E$ in $Q$, and let $b$ be its cut-vertex. Every inner vertex of $E$ has degree at least $k-1$ in $E$. Since $G$ is $2$-connected, some inner vertex $u$ of $E$ has a neighbor on $C$. By \Cref{thm:paths} with exceptional vertex $b$, $E$ has $k-2$ admissible $(u,v)$-paths for every $v\in V(E)\setminus\{u\}$. For $v\in V(Q)\setminus V(E)$, append a fixed $(b,v)$-path in $Q$ avoiding the inner vertices of $E$ to the $k-2$ admissible $(u,b)$-paths in $E$.
\end{proof}

Suppose that $W$ is connected. Choose $u$ as in Claim~\ref{cl:component-paths}, and let $x$ be its neighbor on $C$. Choose $y\in V(C)\setminus\{c\}$ at distance $s-2$ from $x$ on $C$. Since $d_G(y)\ge k$, $y$ has a neighbor $v$ in $W$. Since $u$ has only one neighbor on $C$, $v\ne u$. By Claim~\ref{cl:component-paths}, $W$ has $k-2$ admissible $(u,v)$-paths. The two $(x,y)$-sections of $C$, together with $xu$ and $vy$, form $(u,v)$-paths of lengths $s$ and $s+4$. Combined with the paths in $W$, they form $k$ admissible cycles, a contradiction. Therefore $W$ is disconnected.

Every component of $W$ has at least two neighbors on $C$ since $G$ is $2$-connected. By Claim~\ref{cl:component-paths}, for each component $Q$ of $W$, there are $k-2$ admissible paths between two fixed distinct vertices of $C$, with all internal vertices in $Q$. Choose such path families in two components. Join their endpoint pairs by two fixed vertex-disjoint sections of $C$ whose internal vertices avoid the four endpoints, allowing sections of length zero at common endpoints. Together with these sections, the two path families form $(k-2)+(k-2)-1=2k-5$ admissible cycles, a contradiction if $k\ge5$.
Thus $k=4$. Furthermore, if there is a component $Q$ of $W$ with three admissible paths between two distinct vertices of $C$ with all internal vertices in $Q$, combining them with two admissible paths through another component of $W$ as above, we obtain $3+2-1=4$ admissible cycles, a contradiction. Thus we have the following.\begin{equation}\label{eq:at-most-two}\tag{$\ast$}\parbox{0.88\linewidth}{For every component $Q$ of $W$ and any two distinct vertices of $C$, there are at most two admissible paths between them with all internal vertices in $Q$.}\end{equation} By $4$-weakness, every component $Q$ of $W$ has at least three neighbors on $C$. Otherwise, its two neighbors would form a $2$-cut separating $Q\ne\{c\}$, a contradiction.

\begin{claim}\label{cl:all-pairs}
For every component $Q$ of $W$ and distinct $x,y\in N_G(V(Q))\cap V(C)$, there are two admissible $(x,y)$-paths with all internal vertices in $Q$.
\end{claim}

\begin{proof}
If $Q$ is $2$-connected, then $d_Q(v)\ge3$ for every $v\in V(Q)\setminus\{c\}$, since each vertex of $Q$ has at most one neighbor on $C$. Apply \Cref{thm:paths}, with exceptional vertex $c$, between distinct neighbors of $x,y$ in $Q$ and add the edges to $x,y$.

Suppose that $Q$ is not $2$-connected. By the degree condition, at most one end-block of $Q$ is an edge, namely one with leaf $c$, so $Q$ has a $2$-connected end-block $E_1$. If $Q$ has another $2$-connected end-block, let $E_2$ be one and let $u_2$ be an inner vertex of $E_2$ with a neighbor on $C$. Otherwise, the blocks of $Q$ form a chain with two end blocks $E_1$ and an edge $bc$, and so let $E_2$ be the block other than $bc$ containing $b$ (possibly $E_1=E_2$), which is $2$-connected as $d_G(b)\ge4$, and let $u_2=c$. Let $x_2$ be the neighbor of $u_2$ on $C$. By $4$-weakness, the inner vertices of $E_1$ have at least two neighbors on $C$, so we may choose an inner vertex $u_1$ of $E_1$ with a neighbor $x_1\ne x_2$ on $C$.

Suppose that $E_1\ne E_2$, and let $P$ be a $(u_1,u_2)$-path in $Q$. 
Then each of $E_1$ and $E_2$ contains exactly two vertices of $P$ that are cut-vertices of $Q$ or endpoints of $P$. 
By \Cref{thm:paths}, with exceptional vertex $c$, each of $E_1$ and $E_2$ has two admissible paths between these two vertices. Together with the rest of $P$ and the edges $x_1u_1,u_2x_2$, these form three admissible $(x_1,x_2)$-paths with all internal vertices in $Q$, contrary to \eqref{eq:at-most-two}. Hence $E_1=E_2$, so $Q$ is obtained from $E_1$ by adding the leaf $c$ adjacent to $b$. Write $E=E_1$ and $x_0=x_2$.

If $b$ had a neighbor $y_0$ on $C$, then $x_0cby_0$ and a shortest $(x_0,y_0)$-section of $C$ form a cycle of length at most $s+3\le 2s$. By the minimality of $C$, this is a shortest cycle containing $c$, contrary to the choice of $C$ since $c\in V(Q)$. Thus $b$ has no neighbor on $C$.

For each of $x,y$, choose an edge to $E$, using $x_0cb$ instead if the vertex is $x_0$. By \Cref{thm:paths}, with exceptional vertex $b$, $E$ has two admissible paths between the resulting distinct endpoints in $E$. Together with the paths to $x,y$, they form the required $(x,y)$-paths.
\end{proof}

\begin{claim}\label{cl:two-components}
$W$ has exactly two components.
\end{claim}

\begin{proof}
Let $G'=G$ in Type~I and $G'=G-c+xy$ in Type~II, where $N_G(c)=\{x,y\}$. Then $G'$ is $3$-connected. Let $C'$ be the cycle of $G'$ obtained from $C$ by replacing $xcy$ with $xy$ if $c\in V(C)$ in Type~II, and $C'=C$ otherwise. The components of $G'-V(C')$ correspond to those of $W$, since $c$ has at most one neighbor on $C$ when $c\notin V(C)$. Let $H$ be obtained from $G'$ by contracting each of them to a vertex. By \Cref{lem:contraction}, $H$ is $3$-connected.

Suppose that $W$ has at least three components. Since $H$ is $3$-connected, choose a cycle $D$ through three contracted vertices. In Type~II with $c\in V(C)$, replace $xy$ by $xcy$ if $xy\in E(D)$. At the three selected vertices, use the two admissible paths from Claim~\ref{cl:all-pairs}; at every other contracted vertex on $D$, use a fixed path through its component. Together with the sections of $D$ on $C$, these paths form $2+2+2-2=4$ admissible cycles in $G$, a contradiction.
\end{proof}

By Claim~\ref{cl:two-components}, let $Q_1,Q_2$ be the two components of $W$, and let $A_i=N_G(V(Q_i))\cap V(C)$ for $i=1,2$.
Suppose that there are two vertices $u,v$ in $A_i$ at distance $s-2$ on $C$.
By Claim~\ref{cl:all-pairs}, there are two admissible $(u,v)$-paths with all internal vertices in $Q_i$.
Combining each of these paths with each of the two sections of $C$ between $u$ and $v$ forms four admissible cycles, a contradiction.
Hence, no $A_i$ contains two vertices at distance $s-2$ on $C$ and this implies $|A_i| \leq s$. 
Every vertex of $C$ different from $c$ has a neighbor in $W$. If a vertex of $C$ belongs to both $A_1$ and $A_2$, neither vertex at distance $s-2$ from it has a neighbor in $W$, although at least one is different from $c$, a contradiction. Thus $A_1$ and $A_2$ are disjoint, and $A_1\cup A_2\supseteq V(C)\setminus\{c\}$.
Since $|A_i|\le s$ and $|A_1|+|A_2|\ge 2s-1$, we have $s-1\le|A_i|\le s$.
Thus $C[A_i]$ is disconnected, since otherwise $A_i$ would have two vertices at distance $s-2$. As at most one vertex of $C$ lies outside $A_1\cup A_2$, there are at least three edges of $C$ joining $A_1$ to $A_2$.

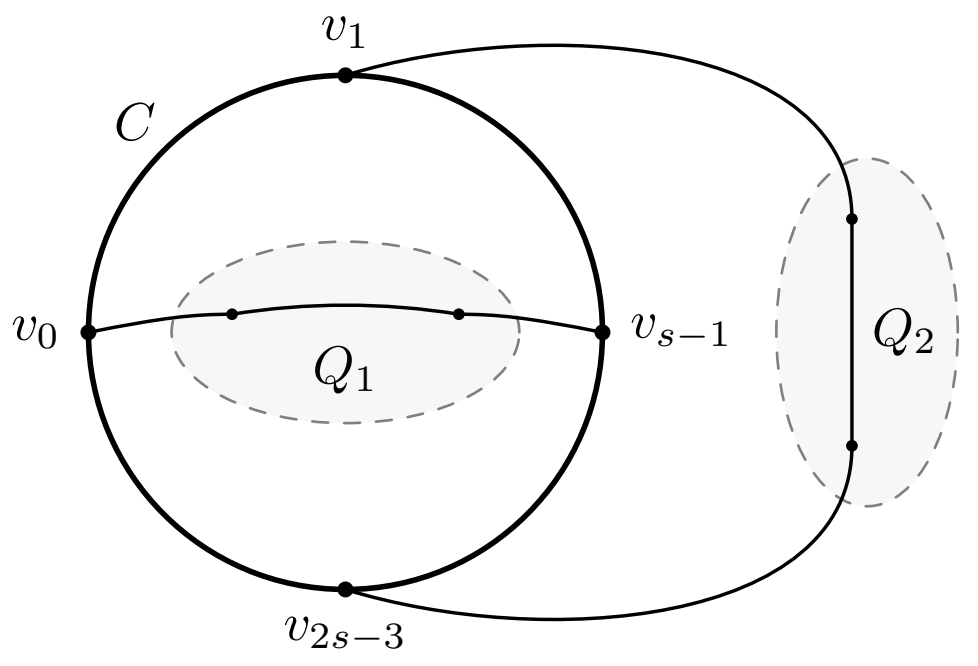
\begin{figure}[htbp]
\centering
\begin{tikzpicture}[font=\small, line cap=round, line join=round]
\coordinate (a) at (-1.7,0);
\coordinate (b) at (0,1.7);
\coordinate (c) at (1.7,0);
\coordinate (d) at (0,-1.7);
\draw[dashed,gray,fill=black!3,line width=.5pt] (0,0) ellipse (1.15 and .6);
\draw[dashed,gray,fill=black!3,line width=.5pt] (3.45,0) ellipse (.6 and 1.15);
\draw[line width=1.05pt] (0,0) circle (1.7);
\node at (135:1.97) {$C$};
\draw[line width=.65pt] (a) .. controls (-1.2,.1) and (-1,.12) .. (-.75,.12) .. controls (-.25,.2) and (.25,.2) .. (.75,.12) .. controls (1,.12) and (1.2,.1) .. (c);
\draw[line width=.65pt] (b) .. controls (1.3,2.1) and (3.35,1.95) .. (3.35,.75) -- (3.35,-.75) .. controls (3.35,-1.95) and (1.3,-2.1) .. (d);
\node at (0,-.25) {$Q_1$};
\node at (3.7,0) {$Q_2$};
\foreach \point in {(-.75,.12),(.75,.12),(3.35,.75),(3.35,-.75)} {\fill \point circle (1.1pt);}
\foreach \vertex in {a,b,c,d} {\fill (\vertex) circle (1.5pt);}
\node[left=4pt,inner sep=1pt] at (a) {$v_0$};
\node[above=4pt,inner sep=1pt] at (b) {$v_1$};
\node[right=4pt,inner sep=1pt] at (c) {$v_{s-1}$};
\node[below=4pt,inner sep=1pt] at (d) {$v_{2s-3}$};
\end{tikzpicture}
\caption{The cycle $C$ and paths through $Q_1,Q_2$ joining the marked vertices in pairs.}\label{fig:cycle-sections}
\end{figure}

Write $C=v_0v_1\cdots v_{2s-1}v_0$. By cyclically relabeling the vertices, we may assume that $v_0v_1$ joins $A_1$ to $A_2$ and that $v_{s-1},v_{2s-3}\ne c$, since at most two of these edges are excluded by the latter conditions. By symmetry, assume $v_0\in A_1$ and $v_1\in A_2$. Since $v_{s-1},v_{2s-3}\in A_1\cup A_2$ and neither $A_i$ contains two vertices at distance $s-2$, we have $v_{s-1}\in A_1$ and $v_{2s-3}\in A_2$.

By Claim~\ref{cl:all-pairs}, there are two admissible $(v_0,v_{s-1})$-paths with all internal vertices in $Q_1$ and two admissible $(v_1,v_{2s-3})$-paths with all internal vertices in $Q_2$. Their length sums form a three-term arithmetic progression with common difference two. Together with the two pairs of sections of $C$ shown in Figure~\ref{fig:cycle-sections}, of total lengths $1+(s-2)=s-1$ and $(s-2)+3=s+1$, these form four admissible cycles, a contradiction.
\end{proof}

\section{Proofs of the main theorems}\label{sec:main-proofs}

\begin{proposition}\label{prop:rooted}
Let $k\ge3$, let $G$ be a $2$-connected bipartite graph, and let $c\in V(G)$ satisfy $d_G(v)\ge k$ for every $v\in V(G)\setminus\{c\}$. Then $G$ has an $(r\bmod 2k)$-cycle for every even integer $r\not\equiv2\pmod{2k}$. Moreover, if $k\ge4$, then either $G$ contains $k$ admissible cycles or $G\in\cE_k$. In particular, for $k\ge4$, $G$ has no $(2\bmod 2k)$-cycle if and only if $G\in\cE_k$.
\end{proposition}

\begin{proof}
Suppose $k=3$. Let $H$ be obtained from two disjoint copies of $G$ by joining the copies of $c$ by an edge. Then $H$ is bipartite and $\delta(H)\ge3$. Since $H$ has no Petersen graph as an induced subgraph,   \Cref{thm:LY} implies that $H$ has a $(4\bmod 6)$-cycle. The added edge is a bridge, so $G$ has a $(4\bmod 6)$-cycle. By \Cref{thm:CS}, $G$ has a $(0\bmod 3)$-cycle, which is a $(0\bmod 6)$-cycle.

Suppose $k\ge4$. If $G$ has $k$ admissible cycles, their lengths are consecutive even integers and represent every even residue modulo $2k$. Otherwise, $G$ is $k$-weak by \Cref{lem:weak}. By \Cref{lem:four-cycle,lem:no-four-cycle}, $G\in\cE_k$. In this case, the cycle lengths are exactly $4,6,\ldots,2k$, so $G$ has an $(r\bmod 2k)$-cycle for every even $r\not\equiv2\pmod{2k}$ and no $(2\bmod 2k)$-cycle.
\end{proof}

\begin{proof}[Proof of \Cref{thm:main}]
Let $B$ be an end-block of $G$. Let $c$ be its cut-vertex, if one exists; otherwise choose any $c\in V(B)$. By the degree condition, $B$ is $2$-connected and $d_B(v)\ge k$ for every $v\in V(B)\setminus\{c\}$. By \Cref{prop:rooted}, $B$ has an $(r\bmod 2k)$-cycle for every even $r\not\equiv2\pmod{2k}$, and has no $(2\bmod 2k)$-cycle if and only if $B\in\cE_k$. Since every cycle is contained in a block, $G$ has no $(2\bmod 2k)$-cycle if and only if every end-block belongs to $\cE_k$ and every non-end-block has no such cycle.
\end{proof}

\begin{proof}[Proof of \Cref{thm:general}]
If $k=2$, the result follows from the theorem of Dean, Lesniak, and Saito~\cite{DLS} that every graph of minimum degree at least three has a $(0\bmod 4)$-cycle.

Suppose $k\ge3$. Choose a bipartition of $V(G)$ maximizing the number of edges between its parts, and let $H$ be the spanning subgraph formed by these edges. By maximality, $d_H(v)\ge\lceil d_G(v)/2\rceil\ge k$ for every $v\in V(G)$. Let $B$ be an end-block of $H$ and let $c$ be its cut-vertex, if one exists; otherwise choose any $c\in V(B)$. Since $\delta(H)\ge k$, $B$ is $2$-connected. Every vertex of $V(B)\setminus\{c\}$ has all its neighbors in $B$, so $d_B(v)\ge k$ for every $v\in V(B)\setminus\{c\}$. By \Cref{prop:rooted}, $B$, and hence $G$, has an $(r\bmod 2k)$-cycle for every even $r\not\equiv2\pmod{2k}$.
\end{proof}

\bigskip

\noindent\textbf{Declaration of AI assistance.} The authors formulated the research problem and developed the mathematical approach. GPT-5.6 Sol (OpenAI) assisted with the exposition under the authors’ guidance. The authors reviewed and revised the manuscript and take full responsibility for its content.

\end{document}